\documentclass[12pt,a4paper,reqno]{amsart}
\usepackage{amsmath,amssymb,amsthm}
\usepackage{enumerate}
\usepackage{url}
\usepackage{fullpage}
\usepackage[nobysame]{amsrefs}

\newcommand{\eps}{\varepsilon}
\newcommand{\D}{\mathfrak{D}}
\newcommand{\M}{\mathfrak{M}}
\newcommand{\m}{\mathfrak{m}}
\newcommand{\tr}{\mathfrak{t}}
\newcommand{\X}{\mathfrak{X}}

\newcommand{\dx}{\mathrm{d}}
\newcommand{\R}{\mathbb{R}}
\newcommand{\Odip}[2]{\mathcal{O}_{#1}\!\left(#2\right)\mathchoice{\!}{}{}{}}
\newcommand{\Odi}[1]{\Odip{}{#1}}
\newcommand{\odip}[2]{{o}_{#1}\!\left(#2\right)\mathchoice{\!}{}{}{}}
\newcommand{\odi}[1]{\odip{}{#1}}

\renewcommand{\qedsymbol}{$\square$}
\newenvironment{Proof}[1][Proof]{\par\noindent\textbf{#1.}~}
{\hfill\qedsymbol\smallskip\par}
\newtheorem{Theorem}{Theorem}
\newtheorem{Lemma}{Lemma}

\allowdisplaybreaks

\title{A Diophantine problem with prime variables}

\author{Alessandro Languasco \& Alessandro Zaccagnini}

\begin{document}


\begin{abstract}
We study the distribution of the values of the form
$\lambda_1 p_1 + \lambda_2 p_2 + \lambda_3 p_3^k$, where $\lambda_1$,
$\lambda_2$ and $\lambda_3$ are non-zero real number not all of the
same sign, with $\lambda_1 / \lambda_2$ irrational, and $p_1$, $p_2$
and $p_3$ are prime numbers.
We prove that, when $1 < k < 4 / 3$, these value approximate rather
closely any prescribed real number.
\end{abstract}

\maketitle

\leftline{Dedicated to Prof.~R.~Balasubramanian on the occasion of his
60th birthday}

\section{Introduction}

The problem that we want to study in this paper can be stated in
general as follows: given $r$ non-zero real numbers $\lambda_1$,
\dots, $\lambda_r$, and positive real numbers $k_1$, \dots, $k_r$,
approximate a given real number $\varpi$ by means of values of the
form
\begin{equation}
\label{values}
  \lambda_1 p_1^{k_1}
  +
  \cdots
  +
  \lambda_r p_r^{k_r},
\end{equation}
where $p_1$, \dots, $p_r$ denote primes.
If $\rho = 1 / k_1 + 1 / k_2 + \dots + 1 / k_r$ is ``small'' the goal
is to show that
\begin{equation}
\label{general-ineq}
  \bigl|
    \lambda_1 p_1^{k_1} + \lambda_2 p_2^{k_2} + \cdots + \lambda_r p_r^{k_r}
    -
    \varpi
  \bigr|
  <
  \eta
\end{equation}
has infinitely many solutions for every fixed $\eta > 0$.
If $\rho$ is ``large'' one expects to be able to prove the stronger
result that, in fact, some $\eta \to 0$ is admissible in
\eqref{general-ineq}: more precisely, it should be possible to take
$\eta$ as a small negative power of $\max_j p_j$.
The number of variables $r$ also plays a role, of course.
Some hypothesis on the irrationality of at least one ratio
$\lambda_i / \lambda_j$ is necessary, and also on signs, if one wants
to approximate to all real numbers and not only some proper subset.
We will make everything precise in due course.

Many such results are known, with various types of assumptions and
conclusions, and we now give a brief description of a few among them.
Vaughan \cite{Vaughan1974a} has $r = 3$ and $k_j = 1$ for all $j$, the
non-zero coefficients $\lambda_j$ not all of the same sign with
$\lambda_1 / \lambda_2$ irrational.
In this case $\eta$ is essentially $(\max_j p_j)^{-1 / 10}$.
The paper \cite{Vaughan1974b} contains more elaborate results of the
same kind, with the same integral exponent $k \ge 1$ for all
primes.
Baker and Harman \cite{BakerH1982} and Harman \cite{Harman1991} have a
result similar to Vaughan's \cite{Vaughan1974a} with
$\eta = (\max_j p_j)^{-1 / 6}$ and $\eta = (\max_j p_j)^{-1 / 5}$
respectively.
This has been recently improved to $\eta = (\max_j p_j)^{-2 / 9}$ by
Matom\"aki \cite{Matomaki2010b}.

The papers by Br\"udern, Cook and Perelli \cite{BrudernCP1997},
Cook \cite{Cook2001}, Cook and Fox \cite{CookF2001},
Harman \cite{Harman2004}, Cook and Harman \cite{CookH2006} all deal
with the number of ``exceptional'' real numbers $\varpi$ that can not
be well approximated by values of type \eqref{values}, but in this
case there are many differences with the results quoted above.
First, $\eta$ does not depend on the primes $p_j$ but rather on
$\varpi$ (it is a small negative power of $\varpi$, in fact), but, in
their setting, this is essentially equivalent to the alternative
statement as we shall see presently.
It is more important, of course, to define carefully what
``exceptional'' means.
Actually, the results apply to suitable sequences of positive real
numbers $\varpi_n$ with limit $+\infty$, and it is shown that the
number of exceptional elements in the sequence, that is, elements that
can not be approximated within the prescribed precision, is small in a
strong quantitative sense.
The assumption is that the coefficients $\lambda_j$ are all positive,
which is not a restriction in this case, that the ratio
$\lambda_1 / \lambda_2$ is irrational and algebraic, and that $k_j$ is
the same positive integer $k$ for all $j$.
The assumption on $\lambda_1 / \lambda_2$ is needed to deal with some
exponential sums on the so-called ``minor'' arc.

Tolev \cite{Tolev1992} has $r = 3$, the coefficients $\lambda_j$ all
equal to $1$ and all the exponents $k_j$ equal to a constant
$k \in (1, 15 / 14)$.
The conclusion is that all sufficiently large real numbers $\varpi$
can be approximated, with $\eta$ a negative power of $\varpi$.

Parsell \cite{Parsell2003} considers two primes and a large number of
powers of $2$, so that in a sense $\rho = 2 + \eps$, but $r$ is large
and $\eta$ is arbitrary but fixed.
This has been improved in \cite{LanguascoZaccagnini2010c} by the
present authors, who showed that a smaller number of powers of $2$ is
needed.
In a similar vein, Languasco and Settimi \cite{LanguascoS2012} have
the corresponding result with one prime, two squares of primes and a
large number of powers of $2$.
Finally, the present authors \cite{LanguascoZ2012c} have a result with
one prime and three squares of primes, so that $\rho = 5 / 2$, $r = 4$
and $\eta = (\max_j p_j)^{- 1 / 18 + \eps}$, while in \cite{LanguascoZ2012d}
they deal with one prime, the square of a prime and the $k$-th power
of a prime with $k \in (1, 33 / 29)$ and
$\eta = ( \max_j p_j )^{-(33 - 29 k) / (72 k) + \eps}$.
In all of these papers, it is assumed that one, carefully chosen,
among the ratios $\lambda_i / \lambda_j$ is irrational.

Our main result is the following Theorem.

\begin{Theorem}
\label{Th:appl}
Let $1 < k < 4 / 3$ be a real number and assume that $\lambda_1$,
$\lambda_2$, and $\lambda_3$ are non-zero real numbers, not all of the
same sign and that $\lambda_1 / \lambda_2$ is irrational.
Let $\varpi$ be any real number.
For any $\eps > 0$ the inequality
\begin{equation}
\label{main-ineq}
  \bigl\vert
    \lambda_1 p_1 + \lambda_2 p_2 + \lambda_3 p_3^k
    -
    \varpi
  \bigr\vert
  \le
  \bigl( \max_j p_j \bigr)^{3/10 - 2/(5k) + \eps}
\end{equation}
has infinitely many solutions in prime variables $p_1$, $p_2$ and
$p_3$.
\end{Theorem}

In the notation above, we have $r = 3$, $\rho = 2 + 1 / k$ and
$\eta = ( \max_j p_j )^{3 / 10 - 2 / (5 k) + \eps}$.
We use the variant of the circle method introduced by Davenport and
Heilbronn \cite{DavenportH1946} to deal with these problems, where the
variables are not necessarily integral.
The following lemmas are the two key ingredients of the proof.
They relate a suitable $L^2$-average of the error on the ``major'' arc
to a generalized version of the so-called Selberg integral, which is a
well-known and widely used tool in this context: see
\cite{BrudernCP1997}, \cite{LanguascoZaccagnini2010c},
\cite{LanguascoS2012}, \cite{LanguascoZ2012c}.
The same argument, with comparatively minor changes, can be used with
$\lambda_1 p_1 + \lambda_2 p_2 + \lambda_3 p_3^2$ (with the same
hypothesis on the ratio $\lambda_1 / \lambda_2$ and on signs as
above), and $\eta = (\max_j p_j)^{-1 / 18 + \eps}$.

Before the statement, we need to define the relevant quantities,
beginning with the exponential sums.
As usual, we write $e(\alpha) = e^{2 \pi i \alpha}$.
For any real $k \ge 1$ we let
\begin{equation}
\label{def-S-U}
  S_k(\alpha)
  =
  \sum_{\delta X \le p^k \le X} \log p \ e(p^k \alpha)
  \qquad\text{and}\qquad
  U_k(\alpha)
  =
  \sum_{\delta X \le n^k \le X} e(n^k \alpha)
\end{equation}
where $\delta$ is a small, fixed positive constant, which may depend
on the coefficients $\lambda_j$.
Then we set
\begin{equation}
\label{k-root-Selberg-int-def}
  J_k(X, h)
  =
  \int_{X}^{2 X}
    \Bigl( \theta((x + h)^{1/k}) - \theta(x^{1/k}) - ((x + h)^{1/k} - x^{1/k})
    \Bigr)^2 \dx x.
\end{equation}
This is the generalized version of the Selberg integral referred to
above: the classical function is $J_1(X, h)$.

\begin{Lemma}
\label{k-root-BCP-Gallagher}
Let $k \geq 1$ be a real number.
For $0 < Y \leq 1/2$  we have
\[
  \int_{-Y}^Y
    \vert S_{k}(\alpha) - U_{k}(\alpha) \vert^2 \dx \alpha
  \ll_k
  \frac{X^{2/k-2}\log^{2}X}{Y}
  +
  Y^2 X
  +
  Y^2 J_{k} \Bigl( X, \frac{1}{2Y} \Bigr),
\]
where $J_{k}(X, h)$ is defined in \eqref{k-root-Selberg-int-def}.
\end{Lemma}

This is Theorem~1 of \cite{LanguascoZ2012d}.

\begin{Lemma}
\label{k-root-Saffari-Vaughan}
Let $k \geq 1$ be a real number and $\eps$ be an arbitrarily small
positive constant.
There exists a positive constant $c_1 = c_{1}(\eps)$, which does not
depend on $k$, such that
\[
  J_{k}(X, h)
  \ll_{k}
  h^2
  X^{2/k - 1}
  \exp \Big( - c_{1} \Big( \frac{\log X}{\log \log X} \Big)^{1/3} \Big)
\]
uniformly for  $X^{1 - 5 / (6 k) + \eps} \leq h \leq X$.
\end{Lemma}

This is the special case $C = 12 / 5$ of Theorem~2 of
\cite{LanguascoZ2012d}.

\section{Proof of Theorem \ref{Th:appl}}

In order to prove that \eqref{main-ineq} has infinitely many
solutions, it is sufficient to show the existence of an increasing
sequence $X_n$ with limit $+\infty$ such that \eqref{main-ineq} has at
least a solution with $\max_j p_j \in [\delta X_n, X_n]$.
This sequence actually depends on rational approximations for
$\lambda_1 / \lambda_2$: more precisely, we recall that there are
infinitely many pairs of integers $a$ and $q$ such that $(a, q) = 1$,
$q > 0$ and
\[
  \Bigl\vert \frac{\lambda_1}{\lambda_2} - \frac aq \Bigr\vert
  \le
  \frac 1{q^2}.
\]
We take the sequence $X = q^{5k / (k+2)}$ (dropping the useless suffix $n$)
and then, as customary, define all of the circle-method parameters in
terms of $X$.
We may obviously assume that $q$ is sufficiently large.
The choice of the exponent $5k / (k+2)$ is justified in the discussion
following the proof of Lemma~\ref{Lemma-approx}.
As usual, we approximate to $S_k$ using the function
\[
  T_k(\alpha)
  =
  \int_{(\delta X)^{1/k}}^{X^{1/k}} e(t^k \alpha) \, \dx t
\]
and notice the simple inequality
\begin{equation}
\label{bd-Ti}
  T_k(\alpha)
  \ll_{k, \delta}
  X^{1 / k -1}
  \min\bigl( X, \vert \alpha \vert^{-1} \bigr).
\end{equation}

Since the variables are not integers, we cannot count exact hits as in
the standard applications of the circle method, only near misses, so
that we need some measure of proximity.
For $\eta > 0$, we detect solutions of \eqref{main-ineq} by means of
the function
\[
  \widehat{K}_{\eta}(\alpha)
  =
  \max(0, \eta - \vert \alpha \vert),
\]
which, as the notation suggests, is the Fourier transform of
\[
  K_{\eta}(\alpha)
  =
  \Bigl( \frac{\sin(\pi \eta \alpha)}{\pi \alpha} \Bigr)^2
\]
for $\alpha \ne 0$, and, by continuity, $K_{\eta}(0) = \eta^2$.
This relation transforms the problem of counting solutions of the
inequality \eqref{main-ineq} into estimating suitable integrals.
We recall the trivial, but crucial, property
\begin{equation}
\label{bd-K(eta)}
  K_{\eta}(\alpha)
  \ll
  \min \Bigl( \eta^2, \vert \alpha \vert^{-2} \Bigr).
\end{equation}

When $\X$ is an interval, a half line, or the union of two such sets
we let
\[
  I(\eta, \varpi, \X)
  =
  \int_{\X} S_1(\lambda_1 \alpha) S_1(\lambda_2 \alpha)
    S_k(\lambda_3 \alpha) K_{\eta}(\alpha)
    e(-\varpi \alpha) \, \dx \alpha.
\]
The starting point of the method is the observation that
\begin{align*}
  I(\eta, \varpi, \R)
  &=
  \sum_{\substack{\delta X \le p_1, p_2 \le X \\ \delta X \le p_3^k \le X}}
    \log p_1 \log p_2 \log p_3
  \int_{\R}
    K_{\eta}(\alpha)
    e \bigl(
        (\lambda_1 p_1 + \lambda_2 p_2 + \lambda_3 p_3^k
        - \varpi) \alpha
      \bigr) \, \dx \alpha \\
  &=
  \sum_{\substack{\delta X \le p_1, p_2 \le X \\ \delta X \le p_3^k \le X}}
    \log p_1 \log p_2 \log p_3
  \max(0,
       \eta - \vert \lambda_1 p_1 + \lambda_2 p_2 + \lambda_3 p_3^k - \varpi
              \vert) \\
  &\le
  \eta (\log X)^3 \mathcal{N}(X),
\end{align*}
where $\mathcal{N}(X)$ denotes the number of solutions of the
inequality \eqref{main-ineq} with $p_1$, $p_2 \in [\delta X, X]$ and
$p_3^k \in [\delta X, X]$.
In other words, $I(\eta, \varpi, \R)$ provides a lower bound for the
quantity that we are interested in.

We now give the definitions that we need to set up the method.
More definitions will be given at appropriate places later.
We let $P = P(X) = X^{5 /(6k) - \eps} $,
$\eta = \eta(X) = X^{3/10 - 2/(5k) + \eps}$, and
$R = R(X) = \eta^{-2} (\log X)^{3/2}$.
The choice for $P$ is justified at the end of \S\ref{subs:J4}, the one
for $\eta$ at the end of \S\ref{sec:intermediate} and the one for $R$
at the end of \S\ref{sec:trivial}.
See also \S\ref{sec:param} for a fuller discussion.
We now decompose $\R$ as $\M \cup \m \cup \tr$ where
\[
  \M
  =
  \Bigl[ -\frac PX, \frac PX \Bigr],
  \qquad
  \m
  =
  \Bigl( -R, -\frac PX \Bigr) \cup \Bigl(\frac PX, R \Bigr),
  \qquad
  \tr
  =
  \R \setminus( \M \cup \m),
\]
so that
\[
  I(\eta, \varpi, \R)
  =
  I(\eta, \varpi, \M)
  +
  I(\eta, \varpi, \m)
  +
  I(\eta, \varpi, \tr).
\]
The sets $\M$, $\m$ and $\tr$ are called the major arc, the
intermediate (or minor) arc and the trivial arc respectively.
In \S\ref{sec:major} we prove that the major arc yields the main term
for $I(\eta, \varpi, \R)$.
We show in \S\ref{sec:intermediate} that the contribution of the
intermediate arc does not cancel the main term, exploiting the
hypothesis that $\lambda_1 / \lambda_2$ is irrational to prove that
$\vert S_1(\lambda_1 \alpha) \vert$ and
$\vert S_1(\lambda_2 \alpha) \vert$ can not both be large for
$\alpha \in \m$: see Lemma~\ref{Lemma-approx} for the details.
The trivial arc, treated in \S\ref{sec:trivial}, only gives a rather
small contribution.

From now on, implicit constants may depend on the coefficients
$\lambda_j$, on $\delta$, $k$ and $\varpi$.

\section{The major arc}
\label{sec:major}

We write
\begin{align*}
  I(\eta, \varpi, \M)
  &=
  \int_{\M}
    S_1(\lambda_1 \alpha) S_1(\lambda_2 \alpha) S_k(\lambda_3 \alpha)
    K_{\eta}(\alpha) e(-\varpi \alpha) \, \dx \alpha \\
  &=
  \int_{\M}
    T_1(\lambda_1 \alpha) T_1(\lambda_2 \alpha) T_k(\lambda_3 \alpha)
    K_{\eta}(\alpha) e(-\varpi \alpha) \, \dx \alpha \\
  &\qquad+
  \int_{\M}
    \bigl( S_1(\lambda_1 \alpha) - T_1(\lambda_1 \alpha) \bigr)
    T_1(\lambda_2 \alpha) T_k(\lambda_3 \alpha)
    K_{\eta}(\alpha) e(-\varpi \alpha) \, \dx \alpha \\
  &\qquad+
  \int_{\M}
    S_1(\lambda_1 \alpha)
    \bigl(S_1(\lambda_2 \alpha) - T_1(\lambda_2 \alpha) \bigr)
    T_k(\lambda_3 \alpha)
    K_{\eta}(\alpha) e(-\varpi \alpha) \, \dx \alpha \\
  &\qquad+
  \int_{\M}
    S_1(\lambda_1 \alpha) S_1(\lambda_2 \alpha)
    \bigl(S_k(\lambda_3 \alpha) - T_k(\lambda_3 \alpha) \bigr)
    K_{\eta}(\alpha) e(-\varpi \alpha) \, \dx \alpha \\
  &=
  J_1 + J_2 + J_3 + J_4,
\end{align*}
say.
We will give a lower bound for $J_1$ and upper bounds for $J_2$, $J_3$
and $J_4$.
For brevity, since the computations for $J_3$ are similar to, but
simpler than, the corresponding ones for $J_2$ and $J_4$, we will skip
them.

\section{Lower bound for $J_1$}
The lower bound $J_1 \gg \eta^2 X^{1+1 / k}$ is  proved in a classical way.
We have
\begin{align*}
  J_1
  &=
  \int_{\M}
    T_1(\lambda_1 \alpha) T_1(\lambda_2 \alpha) T_k(\lambda_3 \alpha)
    K_{\eta}(\alpha) e(-\varpi \alpha) \, \dx \alpha \\
  &=
  \int_{\R}
    T_1(\lambda_1 \alpha) T_1(\lambda_2 \alpha) T_k(\lambda_3 \alpha)
     K_{\eta}(\alpha) e(-\varpi \alpha) \, \dx \alpha \\
  &\qquad+
  \Odi{
  \int_{P / X}^{+\infty}
    \vert
      T_1(\lambda_1 \alpha) T_1(\lambda_2 \alpha) T_k(\lambda_3 \alpha)
    \vert K_{\eta}(\alpha) \, \dx \alpha}.
\end{align*}
Using inequalities \eqref{bd-Ti} and \eqref{bd-K(eta)}, we see that
the error term is
\[
  \ll
  \eta^2
  X^{1/ k -1}
  \int_{P / X}^{+\infty}
    \frac{\dx \alpha}{\alpha^3}
  \ll
  \eta^2 X^{1+1/k} P^{- 2}
  =
  \odi{\eta^2 X^{1 + 1 / k}}.
\]
For brevity, we set $\D = [\delta X, X]^2 \times [(\delta X)^{1/k}, X^{1/k}]$
and rewrite the main term in the form
\begin{align*}
  &
  \idotsint_{\D}
    \int_{\R}
      e \bigl(
      (\lambda_1 t_1 + \lambda_2 t_2 + \lambda_3 t_3^k
       - \varpi) \alpha
        \bigr) \, K_{\eta}(\alpha) \, \dx \alpha
        \, \dx t_1 \, \dx t_2 \, \dx t_3  \\
  &=
  \idotsint_{\D}
    \max(0, \eta
             -
             \vert \lambda_1 t_1 + \lambda_2 t_2 + \lambda_3 t_3^k
                    - \varpi
             \vert)
    \, \dx t_1 \, \dx t_2 \, \dx t_3 .
\end{align*}
We now proceed to show that the last integral is $\gg \eta^2 X^{1+1/k}$.
Apart from trivial permutations or changes of sign, there are
essentially two cases:

\begin{enumerate}

\item
$\lambda_1 > 0$, $\lambda_2 < 0$,  $\lambda_3 < 0$;

\item
$\lambda_1 > 0$, $\lambda_2 > 0$,  $\lambda_3 < 0$.

\end{enumerate}

We briefly deal with the second case, the other one being similar.
A suitable change of variables shows that
\begin{align*}
  J_1
  &\gg
  \idotsint_{\D'}
    \max(0, \eta
            -
            \vert \lambda_1 u_1 + \lambda_2 u_2 + \lambda_3 u_3
            \vert)
    \, \frac{\dx u_1 \, \dx u_2 \, \dx u_3}{u_3^{1-1/k}} \\
  &\gg
  X^{1/k -1}
  \idotsint_{\D'}
    \max(0, \eta
            -
            \vert \lambda_1 u_1 + \lambda_2 u_2 + \lambda_3 u_3
            \vert)
    \, \dx u_1 \, \dx u_2 \, \dx u_3 ,
\end{align*}
where $\D' = [\delta X, (1 - \delta) X]^3$, for large $X$.
For $j = 1$, $2$, let
$a_j = 2 \vert \lambda_3 \vert \delta / \vert \lambda_j \vert$,
$b_j = 3 a_j / 2$ and $\mathfrak{I}_j = [a_j X, b_j X]$.
Notice that if $u_j \in \mathfrak{I}_j$ for $j = 1$, $2$, then
\[
  \lambda_1 u_1 + \lambda_2 u_2
  \in
    \bigl[
    4 \vert \lambda_3 \vert \delta X, 6 \vert \lambda_3 \vert \delta X
  \bigr]
\]
so that, for every such choice of $(u_1, u_2)$, the interval
$[a, b]$ with endpoints
$\pm \eta / \vert \lambda_3 \vert +
(\lambda_1 u_1 + \lambda_2 u_2) / \vert \lambda_3 \vert$
is contained in $[\delta X, (1 - \delta) X]$.
In other words, for $u_3 \in [a, b]$ the values of
$\lambda_1 u_1 + \lambda_2 u_2 + \lambda_3 u_3 $
cover the whole interval $[-\eta, \eta]$.
Hence, for any $(u_1, u_2) \in \mathfrak{I}_1 \times \mathfrak{I}_2$
we have
\[
  \int_{\delta X}^{(1 - \delta) X}
    \max(0, \eta
            -
            \vert \lambda_1 u_1 + \lambda_2 u_2 + \lambda_3 u_3
            \vert)
    \, \dx u_3
  =
  \vert \lambda_3 \vert^{-1}
  \int_{-\eta}^{\eta} \max(0, \eta - \vert u \vert) \, \dx u
  \gg
  \eta^2.
\]
Finally,
\[
  J_1
  \gg
  \eta^2
  X^{1 / k -1}
  \iint_{\mathfrak{I}_1 \times \mathfrak{I}_2}
    \dx u_1 \, \dx u_2
  \gg
  \eta^2 X^{1 +  1/ k},
\]
which is the required lower bound.

\section{Bound for $J_2$}

We recall definition \eqref{def-S-U} and notice that the Euler
summation formula implies that
\begin{equation}
\label{bd-T-U}
  T_k(\alpha)
  -
  U_k(\alpha)
  \ll
  1 + \vert \alpha \vert X
  \qquad\text{for any $k \ge 1$.}
\end{equation}
Using \eqref{bd-K(eta)} we see that
\begin{align*}
  J_2
  &\ll
  \eta^2
  \int_{\M}
    \bigl\vert S_1(\lambda_1 \alpha) - T_1(\lambda_1 \alpha) \bigr\vert \,
    \vert T_1(\lambda_2 \alpha) \vert \,
    \vert T_k(\lambda_3 \alpha) \vert \, \dx \alpha \\
  &\le
  \eta^2
  \int_{\M}
    \bigl\vert S_1(\lambda_1 \alpha) - U_1(\lambda_1 \alpha) \bigr\vert \,
    \vert T_1(\lambda_2 \alpha) \vert \,
    \vert T_k(\lambda_3 \alpha) \vert \, \dx \alpha \\
  &\qquad+
  \eta^2
  \int_{\M}
    \bigl\vert U_1(\lambda_1 \alpha) - T_1(\lambda_1 \alpha) \bigr\vert \,
    \vert T_1(\lambda_2 \alpha) \vert \,
    \vert T_k(\lambda_3 \alpha) \vert \, \dx \alpha \\
  &=
  \eta^2
  (A_2 + B_2),
\end{align*}
say.
In order to estimate $A_2$ we use Lemmas \ref{k-root-BCP-Gallagher}
and \ref{k-root-Saffari-Vaughan}.
By the Cauchy inequality and \eqref{bd-Ti} above, for any fixed
$A > 0$ we have
\begin{align*}
  A_2
  &\ll
  \Bigl(
    \int_{-P / X}^{P / X}
      \bigl\vert S_1(\lambda_1 \alpha) - U_1(\lambda_1 \alpha) \bigr\vert^2
      \, \dx \alpha
  \Bigr)^{1 / 2}
  \Bigl(
    \int_{-P / X}^{P / X}
      \vert T_1(\lambda_2 \alpha) \vert^2 \,
      \vert T_k(\lambda_3 \alpha) \vert^2 \,  \dx \alpha
  \Bigr)^{1 / 2} \\
  &\ll
    \Bigl(\frac X{(\log X)^A} \Bigr)^{1 / 2}
    \Bigl(
    \int_0^{1 / X} X^{2+2/k} \, \dx \alpha
    +
    \int_{1 / X}^{P / X} \frac{X^{2/k-2}}{\alpha^4} \dx \alpha
  \Bigr)^{1 / 2}
  \ll_A
  \frac{X^{1+ 1 / k}}{(\log X)^{A / 2}}
\end{align*}
by Lemma \ref{k-root-Saffari-Vaughan}, which we can use provided that
$X / P \ge X^{1 / 6 + \eps} $, that is, $P \le X^{5 / 6 - \eps}$.
This proves that $\eta^2 A_2 = \odi{\eta^2 X^{1+ 1 / k}}$.
Furthermore, using inequalities \eqref{bd-Ti} and \eqref{bd-T-U} we
see that
\begin{align*}
  B_2
  &\ll
  \int_0^{1 / X}
    \vert T_1(\lambda_2 \alpha) \vert \,
    \vert T_k(\lambda_3 \alpha) \vert \, \dx \alpha
  +
  X
  \int_{1 / X}^{P / X}
    \alpha \,
    \vert T_1(\lambda_2 \alpha) \vert \,
    \vert T_k(\lambda_3 \alpha) \vert \,  \dx \alpha \\
  &\ll
  \frac1X X^{1 + 1/ k}
  +
  X^{1/k}
  \int_{1 / X}^{P / X}  \, \frac{\dx \alpha}{\alpha}
  \ll
  X^{1 / k} \log P,
\end{align*}
so that $\eta^2 B_2 = \odi{\eta^2 X^{1 + 1 /k }}$.

\section{Bound for $J_4$}
\label{subs:J4}

Inequality \eqref{bd-K(eta)} implies that
\begin{align*}
  J_4
  &\ll
  \eta^2
  \int_{\M}
    \bigl\vert S_1(\lambda_1 \alpha) \bigr\vert \,
    \bigl\vert S_1(\lambda_2 \alpha) \bigr\vert \,
    \bigl\vert S_k(\lambda_3 \alpha) - T_k(\lambda_3 \alpha) \bigr\vert
    \, \dx \alpha \\
  &\ll
  \eta^2
  \int_{\M}
    \bigl\vert S_1(\lambda_1 \alpha) \bigr\vert \,
    \bigl\vert S_1(\lambda_2 \alpha) \bigr\vert \,
    \bigl\vert S_k(\lambda_3 \alpha) - U_k(\lambda_3 \alpha) \bigr\vert
    \, \dx \alpha \\
  &\qquad+
  \eta^2
  \int_{\M}
    \bigl\vert S_1(\lambda_1 \alpha) \bigr\vert \,
    \bigl\vert S_1(\lambda_2 \alpha) \bigr\vert \,
    \bigl\vert U_k(\lambda_3 \alpha) - T_k(\lambda_3 \alpha) \bigr\vert
    \, \dx \alpha \\
  &=
  \eta^2 (A_4 + B_4),
\end{align*}
say.
The Parseval inequality and trivial bounds yield, for any fixed
$A > 0$,
\begin{align*}
  A_4
  &\ll
  X
  \Bigl(
    \int_{\M} \bigl\vert S_1(\lambda_1 \alpha) \bigr\vert^2 \, \dx \alpha
  \Bigr)^{1 / 2}
  \Bigl(
    \int_{\M}
      \bigl\vert S_k(\lambda_3 \alpha) - U_k(\lambda_3 \alpha) \bigr\vert^2
      \, \dx \alpha
  \Bigr)^{1 / 2} \\
  &\ll
  X (X \log X)^{1 / 2}
  \frac PX J_k\Bigl(X, \frac XP \Bigr)^{1 / 2}
  \ll_A
  X^{1+1/k } (\log X)^{1 / 2 - A / 2}
\end{align*}
by Lemmas \ref{k-root-BCP-Gallagher} and \ref{k-root-Saffari-Vaughan}
which we can use provided that $X / P \ge X^{1-5 /(6k)  + \eps}$, that is,
$P \le X^{5 / (6k) - \eps}$.
This proves that $\eta^2 A_4= \odi{\eta^2 X^{1 + 1/ k}}$.
Furthermore, using \eqref{bd-T-U}, the Cauchy inequality and trivial
bounds we see that
\begin{align*}
  B_4
  &\ll
  \int_0^{1 / X}
    \bigl\vert S_1(\lambda_1 \alpha) \bigr\vert \,
    \bigl\vert S_1(\lambda_2 \alpha) \bigr\vert \,  \dx \alpha
  +
  X
  \int_{1 / X}^{P / X}
    \alpha
    \bigl\vert S_1(\lambda_1 \alpha) \bigr\vert \,
    \bigl\vert S_1(\lambda_2 \alpha) \bigr\vert \,  \dx \alpha \\
  &\ll
  X
  +
  P
  \Bigl(
    \int_{1 / X}^{P / X}
      \bigl\vert S_1(\lambda_1 \alpha) \bigr\vert^2 \, \dx \alpha
     \int_{1 / X}^{P / X}
      \bigl\vert S_1(\lambda_2 \alpha) \bigr\vert^2 \, \dx \alpha
  \Bigr)^{1 / 2}
  \ll
  P X \log X.
\end{align*}
Hence $B_4 \ll P X \log X$, so that taking
$P = \odi{X^{1 / k} (\log X)^{-1}}$ we get
$\eta^2 B_4 = \odi{\eta^2 X^{1 + 1 / k}}$.
We may therefore choose
\begin{equation}
\label{P-choice}
  P
  =
  X^{5 / (6 k) - \eps}.
\end{equation}

\section{The intermediate arc}
\label{sec:intermediate}

We need to show that $\vert S_1(\lambda_1 \alpha) \vert$ and
$\vert S_1(\lambda_2 \alpha) \vert$ can not both be large for
$\alpha \in \m$, exploiting the fact that $\lambda_1 / \lambda_2$ is
irrational.
We achieve this using a famous result by Vaughan about $S_{1}(\alpha)$.

\begin{Lemma}[Vaughan \cite{Vaughan1997}, Theorem 3.1]
\label{Vaughan-estim}
Let $\alpha$ be a real number and $a,q$ be positive integers
satisfying $(a, q) = 1$ and $\vert \alpha -a/q \vert < q^{-2}$.
Then
\[
  S_{1}(\alpha)
  \ll
  \Bigl(\frac{X}{\sqrt{q}} + \sqrt{Xq} + X^{4/5} \Bigr) \log^4 X.
\]
\end{Lemma}

\begin{Lemma}
\label{Lemma-approx}
Let $1 \le k < 4 / 3$.
Assume that $\lambda_1 / \lambda_2$ is irrational and let $X = q^{5k / (k+2)}$,
where $q$ is the denominator of a convergent of the continued fraction
for $\lambda_1 / \lambda_2$.
Let $V(\alpha) =
\min \bigl( \vert S_1(\lambda_1\alpha) \vert, \vert S_1(\lambda_2\alpha) \vert \bigr)$.
Then we have
\[
  \sup_{\alpha \in \m} V(\alpha)
  \ll
  X^{4 / 5 + 1/(10 k)}  \log ^{4} X.
\]
\end{Lemma}

\begin{Proof}
Let $\alpha \in \m$ and $Q = X^{2/5 -1/ (5k)}  \leq P$.
By Dirichlet's Theorem, there exist integers $a_{i},q_{i}$  with
$1\leq q_{i}\leq X/Q$ and $(a_{i},q_{i})=1$, such that
$\vert \lambda_{i} \alpha q_{i}-a_{i}\vert \leq Q/X$, for $i=1,2$.
We remark that $a_{1}a_{2} \neq 0$ otherwise we would have $\alpha\in \M$.
Now suppose that  $q_{i} \leq Q$ for $i=1,2$. In this case we get
\[
a_{2}q_{1} \frac{\lambda_{1}}{\lambda_{2}} - a_{1}q_{2}
=
( \lambda_{1} \alpha q_{1}-a_{1}) \frac{a_{2}}{\lambda_{2} \alpha}
-
( \lambda_{2} \alpha q_{2}-a_{2}) \frac{a_{1}}{\lambda_{2} \alpha}
\]
and hence
\begin{equation}
\label{bd-1}
  \left\vert
    a_{2}q_{1} \frac{\lambda_{1}}{\lambda_{2}} - a_{1}q_{2}
  \right\vert
  \leq
  2
  \left(
    1+ \left\vert  \frac{\lambda_{1}}{\lambda_{2}} \right\vert
  \right)
  \frac{Q^{2}}{X}
  <
  \frac{1}{2q}
\end{equation}
for sufficiently large $X$.
Then, from the law of best approximation and the definition of $\m$, we obtain
\begin{equation}
\label{bd-2}
  X^{(k+2)/(5k)}=q
  \leq
  \vert a_{2}q_{1} \vert
  \ll q_{1}q_{2} R
  \leq Q^{2} R
%  \leq X^{4 / 5  -2/(5k) - 3/5 +4/(5k)  } \log^{3/2-2} X,
  \leq X^{(k+2)/(5k) -\eps} ,
\end{equation}
which is absurd.
Hence either $q_{1}>Q$ or $q_{2}>Q$.
 Assume that $q_{1}>Q$.
Using Lemma \ref{Vaughan-estim} on $S_1(\lambda_1 \alpha)$, we have
\begin{align}
\notag
V(\alpha)
\leq
\vert
S_1(\lambda_1 \alpha)
\vert
& \ll
\sup_{Q<q_{1}\leq X/Q }
\left(
\frac{X}{\sqrt{q_{1}}}
+
\sqrt{Xq_{1}}
+
X^{4/5}
\right)
\log^{4}X
\\
\notag
%\label{second-minor}
&
\ll
X^{4/5 + 1/(10k)}
(\log X)^{4}.
\end{align}
The other case is totally similar and hence
Lemma \ref{Lemma-approx} follows.
\end{Proof}

\begin{Lemma}
\label{Lemma:bd-minor}
For $j = 1$ and $2$ we have
\[
  \int_{\m}
    \vert S_1(\lambda_j \alpha) \vert^2 K_{\eta}(\alpha) \, \dx \alpha
  \ll
  \eta X \log X
  \quad\text{and}\quad
  \int_{\m}
    \vert S_k(\lambda_3 \alpha) \vert^2 K_{\eta}(\alpha) \, \dx \alpha
  \ll
  \eta X^{1/k} (\log X)^3.
\]
\end{Lemma}

\begin{Proof}
We have to split the range $[P / X, R]$ into two intervals in
order to use \eqref{bd-K(eta)} efficiently.
In the first case we have
\[
  \int_{\m}
    \vert S_1(\lambda_j \alpha) \vert^2 K_{\eta}(\alpha) \, \dx \alpha
  \ll
  \eta^2
  \int_{P / X}^{1 / \eta}
    \vert S_1(\lambda_j \alpha) \vert^2 \, \dx \alpha
  +
  \int_{1 / \eta}^R
    \vert S_1(\lambda_j \alpha) \vert^2 \, \frac{\dx \alpha}{\alpha^2}
\]
by \eqref{bd-K(eta)}, for $j = 1$, $2$.
By periodicity
\[
  \eta^2
  \int_{P / X}^{1 / \eta}
    \vert S_1(\lambda_j \alpha) \vert^2 \, \dx \alpha
  \ll
  \eta
  \int_0^1
    \vert S_1(\alpha) \vert^2 \, \dx \alpha
  \ll
  \eta X \log X,
\]
by the Prime Number Theorem (PNT).
We also have
\[
  \int_{1 / \eta}^R
    \vert S_1(\lambda_j \alpha) \vert^2 \, \frac{\dx \alpha}{\alpha^2}
  \ll
  \int_{|\lambda_j| / \eta}^{+\infty}
    \vert S_1(\alpha) \vert^2 \, \frac{\dx \alpha}{\alpha^2}
  \ll
  \sum_{n \ge |\lambda_j| / \eta}
    \frac1{(n - 1)^2}
    \int_{n - 1}^n
      \vert S_1(\alpha) \vert^2 \, \dx \alpha
  \ll
  \eta X \log X,
\]
again by the PNT.
This proves the first part of the statement.
For the second part, we argue in a similar way, replacing the PNT by
an appeal to (iii) of Lemma~7 in Tolev
\cite{Tolev1992}.
\end{Proof}

Now let
\begin{align*}
  \X_1
  &=
  \{ \alpha \in [P / X, R] \colon
    \vert S_1(\lambda_1 \alpha) \vert
    \le
    \vert S_1(\lambda_2 \alpha) \vert \} \\
  \X_2
  &=
  \{ \alpha \in [P / X, R] \colon
    \vert S_1(\lambda_1 \alpha) \vert
    \ge
    \vert S_1(\lambda_2 \alpha) \vert \}
\end{align*}
so that $[P / X, R] = \X_1 \cup \X_2$ and
\[
  \Bigl\vert I(\eta, \varpi, \m) \Bigr\vert
  \ll
  \Bigl( \int_{\X_1} + \int_{\X_2} \Bigr)
    \bigl\vert
      S_1(\lambda_1 \alpha) S_1(\lambda_2 \alpha)
      S_k(\lambda_3 \alpha)
    \bigr\vert
    K_{\eta}(\alpha) \, \dx \alpha.
\]
Cauchy's inequality gives
\begin{align*}
  \int_{\X_1}
  &\le
  \max_{\alpha \in \X_1} \vert S_1(\lambda_1 \alpha) \vert
  \Bigl(
    \int_{\X_1}
      \vert S_1(\lambda_2 \alpha) \vert^2 K_{\eta}(\alpha) \, \dx \alpha
  \Bigr)^{1 / 2}
  \Bigl(
    \int_{\X_1}
      \vert S_k(\lambda_3 \alpha) \vert^2 K_{\eta}(\alpha) \, \dx \alpha
  \Bigr)^{1 / 2} \\
  &\ll
  X^{4 / 5 + 1/(10k) }
  (\log X)^{4} (\eta X \log X)^{1 / 2} (\eta X^{1/k} (\log X)^3)^{1 / 2} \\
  &\ll
  \eta X^{13 / 10  +3 / (5k)} (\log X)^{6}
\end{align*}
by Lemmas~\ref{Lemma-approx} and \ref{Lemma:bd-minor}.
The computation on $\X_2$ is similar and gives the same final result.
Summing up,
\[
  \Bigl\vert I(\eta, \varpi, \m) \Bigr\vert
  \ll
  \eta X^{13 / 10  +3 / (5k)} (\log X)^{6},
\]
and this is $\odi{\eta^2 X^{1 + 1/k}}$ provided that
\begin{equation}
\label{eta-choice}
  \eta
  =
  \infty(X^{3/10 - 2 / (5k)} (\log X)^6).
\end{equation}

\section{The trivial arc}
\label{sec:trivial}

Using the Cauchy inequality and a trivial bound for
$S_k(\lambda_3 \alpha)$ we see that
\begin{align*}
  \Bigl\vert I(\eta, \varpi, \tr) \Bigr\vert
  &\le
  2
  \int_R^{+\infty}
    \vert S_1(\lambda_1 \alpha) \vert \, \vert S_1(\lambda_2 \alpha) \vert \,
    \vert S_k(\lambda_3 \alpha) \vert \,
    K_{\eta}(\alpha) \, \dx \alpha \\
  &\ll
  X^{1 / k}
  \Bigl(
    \int_R^{+\infty}
      \vert S_1(\lambda_1 \alpha) \vert^2 K_{\eta}(\alpha) \, \dx \alpha
  \Bigr)^{1/2}
  \Bigl(
    \int_R^{+\infty}
    \vert S_1(\lambda_2 \alpha) \vert^2\, K_{\eta}(\alpha) \, \dx \alpha
  \Bigr)^{1/2} \\
  &\ll
  X^{1 / k}
  C_1^{1 / 2}
  C_2^{1 / 2},
\end{align*}
say, where in the last but one line we used the inequality
\eqref{bd-K(eta)}, and, for $j=1,2$, we set
\[
  C_j
  =
  \int_{\vert \lambda_j \vert R}^{+\infty}
    \frac{\vert S_1(\alpha) \vert^2}{\alpha^2} \, \dx \alpha.
\]
We argue as in the proof of Lemma~\ref{Lemma:bd-minor}.
Using the PNT we have
\[
  C_j
  \ll
  \sum_{n \ge \vert \lambda_j \vert R}
    \frac 1{(n - 1)^2}
    \int_{n - 1}^n \vert S_1(\alpha) \vert^2 \, \dx \alpha
  \ll
  \frac{X \log X}{\vert \lambda_j \vert R}.
\]
Collecting these estimates, we conclude that
\begin{equation}
%\label{final-est-trivial}
\notag
  \Bigl\vert I(\eta, \varpi, \tr) \Bigr\vert
  \ll
  \frac{X^{1+ 1/k} \log X}R.
\end{equation}
Hence, the choice
\begin{equation}
\label{R-choice}
  R
  =
  \eta^{-2} (\log X)^{3/2}
\end{equation}
is admissible.

\section{Remark on the choice of the parameters}
\label{sec:param}

The choice $X = q^{5 k / (k + 2)}$ with $1 \le k < 4 / 3$ arises from
the bounds \eqref{bd-1} and \eqref{bd-2}.
Their combination prevents us from choosing the optimal value
$X = q^2$.
This is justified as follows: neglecting log-powers, let
$X = q^{a(k)}$, $Q = X^{b(k)}$, $\eta = X^{-c(k)}$, and recall the
choices $P = X^{5 / (6 k) - \eps}$ in \eqref{P-choice} and
$R = \eta^{-2} (\log X)^{3 / 2}$ in \eqref{R-choice} which are due,
respectively, to the bound for $B_4$ and for the trivial arc.
Then, essentially, we have to maximize $k$ subject to the constraints
\[
  \begin{cases}
    a(k) \ge 1 \\
    0 \le b(k) \le 5 / (6 k) \\
    c(k) \ge 0 \\
    2 b(k) - 1 \le - 1 / a(k)
    &\text{by \eqref{bd-1},} \\
    2 b(k) + 2 c(k) \le 1 / a(k)
    &\text{by \eqref{bd-2},} \\
    - c(k) \ge \frac12  - \frac1{2 k} - \frac12 b(k)
    &\text{by \eqref{eta-choice},}
  \end{cases}
\]
which is a linear optimization problem in the variables $1 / a(k)$,
$b(k)$, $c(k)$ and $1 / k$.
The solution for this problem is $1 / a(k) = (k + 2) / (5 k)$,
$b(k) = (2 k - 1) / (5 k)$, $c(k) = (4 - 3 k) / (10 k)$,
for $1 / k \ge 3 / 4$, and this is equivalent to the statement of the
main Theorem.

%
%\bibliographystyle{amsplain}
%\bibliography{teonum}
% \bib, bibdiv, biblist are defined by the amsrefs package.

\bigskip

\author{Alessandro LANGUASCO\\
Universit\`a di Padova\\
Dipartimento di Matematica\\
Via Trieste 63\\
35121 Padova, Italy\\
E-mail: languasco@math.unipd.it}

\bigskip

\author{Alessandro ZACCAGNINI \\
Universit\`a di Parma \\
Dipartimento di Matematica \\
Parco Area delle Scienze, 53/a \\
Campus Universitario \\
43124 Parma, Italy \\
E-mail: alessandro.zaccagnini@unipr.it}

\end{document}